\pdfoutput=1

\makeatletter
\IfFileExists{./numapde-latex/numapde-packages.sty}{\providecommand*{\input@path}{}\edef\input@path{{./numapde-latex/}\input@path}}{}
\makeatother

\documentclass{numapde-preprint}

\usepackage{numapde-semantic}
\usepackage{numapde-style}
\usepackage{numapde-local}

\IfFileExists{./numapde-bibliography/numapde.bib}{\addbibresource{./numapde-bibliography/numapde.bib}}{\addbibresource{numapde.bib}}

\hypersetup{
	pdftitle={An Optimal Control Problem for Single-Spot Pulsed Laser Welding},
	pdfauthor={Roland Herzog, Dmytro Strelnikov},
	pdfkeywords={optimal control, pulsed laser welding, solidification velocity, quasilinear heat equation, phase transition}
}

\title{An Optimal Control Problem for Single-Spot Pulsed Laser Welding\thanks{This work was funded by the Federal Ministry for Economic Affairs and Energy (BMWi) of Germany through IGF grant 20.826~B (project OptiPuls), which is gratefully acknowledged.}}
\subtitle{}
\shorttitle{Optimal Control of Single-Spot Pulsed Laser Welding}

\author{Roland Herzog\thanks{Interdisciplinary Center for Scientific Computing, Heidelberg University, 69120 Heidelberg, Germany (\email{roland.herzog@iwr.uni-heidelberg.de}, \url{https://scoop.iwr.uni-heidelberg.de}, \orcid{0000-0003-2164-6575}).}
\and
Dmytro Strelnikov\thanks{Technische Universität Chemnitz, Faculty of Mathematics, 09107 Chemnitz, Germany (\email{dmytro.strelnikov@mathematik.tu-chemnitz.de}, \url{https://www.tu-chemnitz.de/mathematik/part_dgl/people/strelnikov}, \orcid{0000-0002-7668-3640}).}}
\shortauthor{R. Herzog and D. Strelnikov}

\dedication{}

\IfFileExists{numapde-uncolorizeHyperrefIfChanges.sty}{\RequirePackage{numapde-uncolorizeHyperrefIfChanges}}{}

\begin{document}
\maketitle

\begin{abstract}
We consider an optimal control problem for a single-spot pulsed laser welding problem.
The distribution of thermal energy is described by a quasilinear heat equation.
Our emphasis is on materials which tend to suffer from hot cracking when welded, such as aluminum alloys.
A simple precursor for the occurrence of hot cracks is the velocity of the solidification front.
We therefore formulate an optimal control problem whose objective contains a term which penalizes excessive solidification velocities.
The control function to be optimized is the laser power over time, subject to pointwise lower and upper bounds.
We describe the finite element discretization of the problem and a projected gradient scheme for its solution.
Numerical experiments for material data representing the EN~AW~6082-T6 aluminum alloy exhibit interesting laser pulse patterns which perform significantly better than standard ramp-down patterns.\end{abstract}


\begin{AMS}
\href{https://mathscinet.ams.org/msc/msc2010.html?t=49K20}{49K20}, \href{https://mathscinet.ams.org/msc/msc2010.html?t=35K59}{35K59}, \href{https://mathscinet.ams.org/msc/msc2010.html?t=49M05}{49M05}
\end{AMS}

\section{Introduction}
\label{sec:introduction}

Pulsed laser welding is a standard technology to merge metal or thermoplastic components.
Its advantages are the narrow spatial concentration and high peak power of the the heat source, as well as the opportunity to quickly and frequently adjust the laser power in time.
However, in comparison to continuous wave laser processes, pulse laser welding is reported to have an elevated tendency to produce hot cracks during the solidification phase due to higher cooling and thus strain rates.
While small hot cracks do not necessarily affect the strength of the welding seam, they may impair the air- and water-tightness.
Avoiding hot cracks is particularly difficult for the welding of certain aluminum alloys, \eg, some of the 2XXX, 5XXX and most of the 6XXX series, which remains a challenging engineering problem \cite{Katayama:2001:1,ZhangWeckmanZhou:2008:1,BieleninBergmann:2017:1}.

Previous analyses have shown the potential to reduce hot cracking by varying the laser power profile in pulsed laser welding; see, \eg, \cite{BieleninBergmann:2017:1,JiaZhangYuShiLiuWuYeWangTian:2021:1}.
In this paper, we propose an optimal control approach to find power profiles which are optimal in a certain sense.
We concentrate on single-spot pulsed laser welding problems with a view towards aluminum alloy welding.
Since welding seams consist of multiple, partially overlapping welding spots, this work constitutes a significant first step towards the optimization of entire welding seams.

In order to obtain a sufficiently realistic forward model of heat distribution, we need to take into account several physical effects, including temperature dependent heat capacity and thermal conductivity, the enthalpy of fusion and convective heat transfer.
From the mathematical point of view, this results in a quasilinear heat equation.
Evaporation of metal will be disregarded, as well as fluidic motion inside the weld pool.
The thermal energy incurred through the laser into the welded component is modeled through a heat flux boundary condition.
Our objective or cost functional takes into account, among other things, the speed of solidification in order to avoid or reduce the appearance of the hot cracks.

The emphasis of our contribution lies with the description of the quasilinear heat equation model, the formulation of an appropriate cost function, as well as the numerical solution of a discretized version of the optimal control problem by a projected gradient descent scheme. One of the terms of the objective functional which penalizes excessive solidification velocities is rather non-standard and was designed specifically for this problem.

The material is structured as follows.
In \cref{sec:modelling}, we discuss the quasilinear heat equation representing the forward model.
The optimal control problem is described in \cref{sec:optimal_control_problem}.
Its discretization is detailed in \cref{sec:discretization}, where we also present a reduction of the three-dimensional setup to the radially symmetric case.
\Cref{sec:numericals} is devoted to the presentation of optimized laser pulse profiles under various conditions.

\section{Modelling}
\label{sec:modelling}

\begin{figure}[ht]
	\centering
	\newcommand{\varR}{3.0}
\newcommand{\varr}{1.0}
\newcommand{\varZ}{2.5}

\tdplotsetmaincoords{70}{0}

\begin{tikzpicture}[tdplot_main_coords]

	\begin{scope}[canvas is xy plane at z=0]
		\draw (\varR, 0) arc [radius=\varR, start angle=0, end angle=-180];
		\draw [very thin, dashed] (\varR, 0) arc [radius=\varR, start angle=0, end angle=180];
	\end{scope}

	\begin{scope}[canvas is xy plane at z=\varZ]
		\draw [fill] node{.} (0,0);
		\draw (0, 0) circle [radius=\varR];
		\draw [very thin, fill=red, fill opacity=0.2] (0, 0) circle [radius=\varr];
	\end{scope}

	\draw ( \varR, 0, 0) -- ( \varR, 0, \varZ);
	\draw (-\varR, 0, 0) -- (-\varR, 0, \varZ);

	\tdplotsetcoord{P1}{ .5*\varR}{90}{-60}
	\tdplotsetcoord{P2}{    \varR}{90}{-60}
	\tdplotsetcoord{P3}{1.4*\varR}{90}{-60}

	\draw [very thin] (0, 0, 0) -- (P2);
	\draw [->, very thin] (P2) -- (P3) node [above] {$r$};
	\draw [very thin] (0, 0, 0) -- (0, 0, \varZ);
	\draw [->, very thin] (0, 0, \varZ) -- (0, 0, 1.6*\varZ) node [left] {$z$};
	\draw [very thin] (0, 0, 0) -- (.6*\varR, 0, 0);
	\tdplotdrawarc[->, very thin]{(0, 0, 0)}{.5*\varR}{-60}{0}{anchor=north west}{$\varphi$}

	\node[left] at (0, 0, \varZ) {$\Gamma_1$};
	\node[left] at (-.5*\varR, 0,  \varZ) {$\Gamma_2$};
	\node[left] at (-.8*\varR, 0, .5*\varZ) {$\Gamma_3$};
	\node[left] at (-.5*\varR, 0,  0) {$\Gamma_4$};

\end{tikzpicture}
	\caption{Cylinder $\Omega$ and its boundaries.}
	\label{fig:cylinder}
\end{figure}
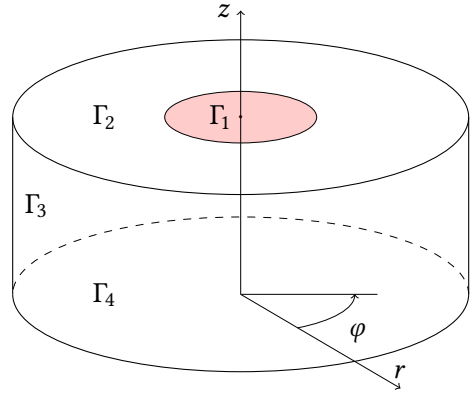

The aim of this section is to derive step-by-step a mathematical model for a single-spot pulsed laser welding problems of aluminum alloys in a cylindrical domain.
To this end, let $\Omega	\subset \R^3$ be an open, orthogonal circular cylinder and $\Gamma = \cup_{i=1}^4 \Gamma_i$ be its boundary surface (see \cref{fig:cylinder}).
Here $\Gamma_1$ is the portion of the boundary affected by the laser beam radiation.
We denote by $\theta(x,t)$ the temperature at the point $x \in \Omega$ at time $t \in [0,T]$.

We are going to describe the temperature evolution inside $\Omega$, and hence the evolution of the welding process, as a solution to a boundary value problem based on the quasilinear heat equation.
The applied nature of the problem provides a few modelling challenges such as temperature dependent properties of the material, liquid/solid phase transition, and a combination of multiple heat transfer mechanisms.
These challenges are sequentially addressed in the following subsections, resulting in a complete model.

\subsection{Enthalpy of Fusion and Volumetric Effective Heat Capacity}
\label{subsec:capacity}

Unlike standard heat dissipation problems when the considered material remains in the same state of matter and its physical properties remain essentially uniform, we deal with a phase transition during the heating and the cooling stages.
These phase transitions are accompanied by an absorption or a release of energy. 
The required amount of additional energy needed to be provided to a specific quantity of the substance to change its state from a solid to a liquid (at constant pressure) is called the \emph{enthalpy of fusion} or the \emph{(latent) heat of fusion}. 
For the opposite transition from a liquid to a solid state \emph{the heat of solidification} has the same absolute value but its sign is reversed.

These phenomena are often modeled in terms of the classical Stefan problem, which is a particular kind of a boundary value problem describing the evolution of a moving boundary between two phases of a material undergoing a phase change; see for instance \cite{Gupta:2003:1}. 
In addition to the underlying heat equation, initial and boundary conditions, the \emph{Stefan condition} is required to provide the energy balance on the phase transition interface.
However, in the present paper we use another approach to integrate the enthalpy of fusion into the boundary value problem. 
Due to the mixed composition of aluminum alloys, we have a wide temperature corridor (rather than a single melting temperature) within which the material melts from a solid to a liquid state. 
The temperature below which the material is fully solid is called \emph{solidus}. 
The temperature above which the material is fully liquid is called \emph{liquidus}. 
In the current study we consider $\text{solidus} = \SI{858}{\K}$ and $\text{liquidus} = \SI{923}{K}$ as reference values.

Considering the above, it becomes more natural in our case to embed the enthalpy of fusion directly into the heat equation by means of \emph{the heat capacity} coefficient. 
In a standard heat dissipation problem with no phase transition, the heat capacity coefficient $c(\theta)$ is a temperature dependent function such that $\int_{\theta_0}^{\theta_1} c(\theta) \, \d\theta$ describes the amount of energy required to heat a unit mass of the material from temperature $\theta_0$ to temperature $\theta_1$. 
In the present model we substitute the heat capacity with an \emph{effective heat capacity} denoted by $c_{\text{eff}}(\theta)$. 
The latter coefficient coincides with $c(\theta)$ outside the solidus--liquidus temperature corridor but has significantly higher values inside, which is meant to achieve the same equality: the total amount of energy required to heat a unit mass of the material from temperature $\theta_0$ to temperature $\theta_1$ (including the enthalpy of fusion if applicable on the interval) is given by the integral $\int_{\theta_0}^{\theta_1} c_{\text{eff}}(\theta) \, \d\theta$.

Another effect to be taken into account is that the \emph{density}~$\rho$ of aluminum alloys changes significantly over the temperature regime under consideration due to thermal expansion and contraction.
However, considering variable volume of the material would lead to a free boundary problem, which significantly increases the complexity of the model.
We therefore take volume changes into account through a temperature dependent density.
Overall, this leads to an \emph{effective volumetric heat capacity} $s(\theta) = c_\textup{eff}(\theta) \rho(\theta)$ in our heat equation.

For the aluminum alloys under consideration, reference values of volumetric heat capacity are given in both the fully solid and the fully liquid state of matter. 
These values show a good linear approximability within a fixed state of matter. 
Therefore, we construct $s(\theta)$ using the following procedure:
\begin{enumerate}
	\item 
		We perform a linear least-squares approximation to the experimental data independently in the solid and in the liquid state of matter.
	\item 
		We choose a $C^1$ cubic spline by filling the liquidus--solidus temperature gap with the uniquely defined cubic polynomial.
	\item 
		In the liquidus--solidus interval, we add without loss of smoothness a cubic spline (consisting of two cubic polynomials) whose integral over the considered interval is equal to the enthalpy of fusion of the selected alloy.
\end{enumerate}

We do not present the above procedure in terms of cumbersome formulas but limit ourselves here to a plot of the resulting effective volumetric heat capacity function; see \cref{fig:coef}.

\subsection{Effective Thermal Conductivity}
\label{subsec:conductivity}

Convective heat transfer in the liquid phase becomes the next modeling challenge caused by the phase transition. 
Due to the Marangoni effect, see \cite{MillsKeeneBrooksShirali:1998:1,Saldi:2012:1}, once the solidus point is passed, the heat transfer in the melting pool significantly increases in radial direction and decreases in axial direction (see \cref{fig:cylinder} for the coordinate axes).
In order to not include the convection term into the core equation we approximate linearly the thermal conductivity coefficient $\kappa(\theta)$ to its measured values in the solid state, and then extrapolate it (separately for the radial and the axial directions) to the temperatures above the liquidus with experimentally selected constants. 
Convective heat transfer in the angular direction is assumed to be zero.
As a result, we have a matrix-valued effective thermal conductivity function $\kappa(\theta) = \diag(\kappa_\textup{ax}(\theta), \kappa_\textup{rad}(\theta),\; 0)$ in diagonal form in a cylindrical coordinate system.

\begin{figure}[ht]
	\centering
	\includegraphics{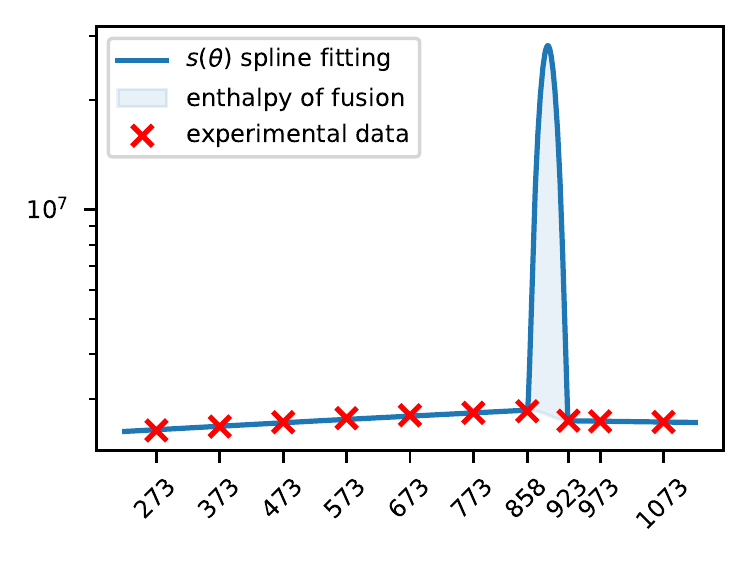}
	\includegraphics{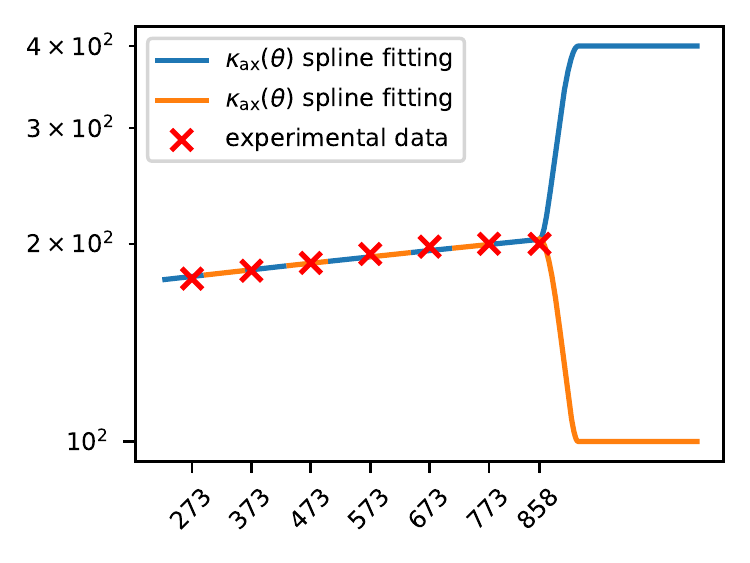}
	\caption{Effective volumetric heat capacity $s(\theta)$ and effective thermal conductivity $\kappa(\theta)$ constructed using a spline fitting procedure.}
	\label{fig:coef}
\end{figure}

The exact algorithm used for constructing $\kappa(\theta)$ (and also $s(\theta)$) can be inspected in code in \cite[\texttt{optipuls.coefficients}]{optipuls_github}. 
We provide a plot of $\kappa_\textup{rad}(\theta)$ and $\kappa_\textup{ax}(\theta)$; see \cref{fig:coef}.
Numerical simulations based on these assumptions have shown reasonable correspondence to the real experiments, \cite{Bielenin:2021:1}.

\subsection{Boundary Conditions}
\label{subsec:boundary_conditions}

While some studies considered the energy introduced by the laser as a volumetric energy source, in this paper we use flux boundary conditions on the boundary part $\Gamma_1$ for this purpose:
\begin{equation*}
	\kappa(\theta(x,t)) \frac{\partial \theta(x,t)}{\partial \vn} 
	= 
	- \eta \, \text{pd}_{\max} \, u(t)
	.
\end{equation*}
Here $\eta$ is the absorption coefficient of the material, $\text{pd}_{\max}$ is the power density of the laser beam, and $u(t)$ is the control function with values in $[0,1]$.
Since the power distribution of the laser beam is taken to be uniform across~$\Gamma_1$ in the current study, the power density is assumed to be constant and hence can be evaluated as the ratio of the maximal total power to the area of the affected spot. 
The absorption coefficient~$\eta$ is also assumed to be constant.

The cooling of the body is a result of the heat flux through the entire boundary except $\Gamma_3$; see \cref{fig:cylinder}.
For simplicity we assume zero heat flux through $\Gamma_3$, which is a reasonable approximation when the radius of $\Omega$ is sufficiently large.
On the remaining parts of the boundary, we distinguish convective and radiative heat fluxes modelled as
\begin{equation*}
	\kappa(\theta(x,t)) \frac{\partial \theta(x,t)}{\partial \vn} 
	= 
	h \cdot (\theta(x,t) - \theta_\textup{amb})
\end{equation*}
and
\begin{equation*}
	\kappa(\theta(x,t)) \frac{\partial \theta(x,t)}{\partial \vn} 
	= 
	k \cdot (\theta(x,t)^4 - \theta^4_\textup{amb})
	,
\end{equation*}
respectively; see for instance \cite[Chapter~3]{Sluzalec:2005:1}.
Here $k = 2.26 \cdot 10^{-9}$ \si{\W\per\m^2\K^4} and $h = 5$ \si{\W\per\m^2} are the radiative and the convective transfer coefficients, respectively, and $\theta_\textup{amb}$ denotes the ambient temperature.

\subsection{Summary of Model Equations}
\label{subsec:equations}

Let us summarize our model based on the above considerations.
We recall that $\Omega\subset \R^3$ is an open, orthogonal circular cylinder and $\Gamma = \cup_{i=1}^4 \Gamma_i$ is its boundary surface as shown in \cref{fig:cylinder}. 
The temperature distribution in $\Omega$ is governed by the quasilinear heat equation
\begin{equation} \label{eq:heat_eq}
	s(\theta(x,t)) \frac{\partial \theta(x,t)}{\partial t} 
	= 
	\div \paren[big](){\kappa(\theta(x,t)) \grad\theta(x,t)}
	,
\end{equation}
where the temperature-dependent coefficients $s(\theta(x,t)) = c_\textup{eff}(\theta(x,t)) \rho(\theta(x,t))$ and $\kappa(\theta(x,t)$ are constructed as $C^1$ cubic splines as detailed in \cite[\texttt{optipuls.coefficients}]{optipuls_github}.

Since we consider single-spot welding, the initial temperature $\theta(x,0)$ inside $\Omega$ is assumed to be constant and equal to the ambience temperature $\theta_\textup{amb}$:
\begin{equation} \label{eq:heat_eq_ic}
	\theta(x,0)
	=
	\theta_\textup{amb}
	\quad
	\text{in } \Omega
	.
\end{equation}
The boundary conditions for \eqref{eq:heat_eq} are
\begin{equation} \label{eq:heat_eq_bc}
	\kappa(\theta(x,t)) \frac{\partial \theta(x,t)}{\partial \vn} 
	= 
	\paren[auto]\{\}{%
		\begin{aligned}
			&
			k (\theta(x,t)^4 - \theta_\textup{amb}^4) + h (\theta(x,t) - \theta_\textup{amb}) - \eta\, \text{pd}_{\max} u(t)
			& 
			& 
			\text{on } \Gamma_1
			, 
			\\
			& 
			k (\theta(x,t)^4 - \theta_\textup{amb}^4) + h (\theta(x,t) - \theta_\textup{amb})
			& 
			& 
			\text{on } \Gamma_2 \cup \Gamma_4
			, 
			\\
			& 
			0
			& 
			& 
			\text{on } \Gamma_3
			.
		\end{aligned}
	}
\end{equation}
We recall that $k$, $h$, $\text{pd}_{\max}$ and $\theta_\textup{amb}$ are known constants.
Moreover, $u(t)$ is the control function with values by $[0,1]$ we seek to determine, which represents the fraction of the maximal laser power to be emitted as a function of time.

\section{Optimal Control Problem}
\label{sec:optimal_control_problem}

This section aims at constructing of an objective functional as a sum of independent penalty terms, each with a different purpose with relation to the single-spot welding application in mind.
As imposed by the application, the desired optimal control representing the emitted laser power profile must:
\begin{enumeratearabic}
	\item 
		provide sufficient welding penetration;
	\item 
		avoid hot cracking during the solidification stage;
	\item 
		ensure complete solidification after welding within the preselected time interval~$[0,T]$;
	\item 
		minimize the total energy consumed by the laser.
\end{enumeratearabic}

In the following subsections we present and discuss four penalty terms designed to target of one these requirements each.
We mention that similar, preliminary ideas were already presented in \cite{BergmannBieleninHerzogHildebrandRiedelSchrickerTrunkWorthmann:2017:1} but with little detail and discussion.

\subsection{Welding Penetration Penalty}
\label{subsec:welding_penetration}

In order to guarantee the \textbf{successful completion of the welding stage} we must ensure that the melting pool has reached a certain predefined depth. 
At the same time, exceeding of this depth would result in an unnecessary increase in energy consumption and the time required for cooling.
Therefore, we select a target point $x_{\text{target}}$ on the symmetry axis of $\Omega$ and a target temperature $\theta_{\text{target}}$ and formulate a term which penalizes the difference between the maximal temperature reached at $x_{\text{target}}$ and the target temperature~$\theta_{\text{target}}$:
\begin{equation} \label{eq:J_penetration}
	J_\textup{penetration} 
	= 
	\frac{\beta_\textup{penetration}}{2} \paren[big](){\norm{\theta(x_{\text{target}},\cdot)}_{L^p(0,T)} - \theta_{\text{target}}}^2
	.
\end{equation}
Here $p$ is sufficiently large so that the $L^p$-norm, which is chosen for simplicity and to avoid non-differentiabilities and state constraints, approximates the $L^\infty$-norm.

\subsection{Solidification Velocity Penalty}
\label{subsec:velocity}

Our main practical goal is to \textbf{avoid the appearance of hot cracks} during the solidification stage. 
As mentioned in \cref{sec:introduction}, we associate hot cracks with high velocities of the solidification front.
We therefore seek to restrict the maximal velocity of the solidification front by introducing a non-standard penalty term derived below.

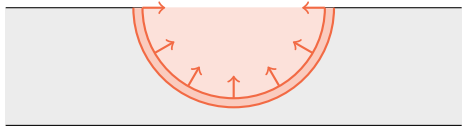
\begin{figure}
	\centering

\newcommand{\varR}{2.5}
\newcommand{\varr}{1}
\newcommand{\vardi}{0.5}
\newcommand{\vardii}{0.8}
\newcommand{\varrlas}{0.4}
\newcommand{\varhlas}{0.8}

\definecolor{mycoral}{RGB}{236,27,75}
\definecolor{myorange}{RGB}{242,106,68}

\begin{tikzpicture}[scale=1.2]

	\draw [draw=none, fill=gray!15] (-\varR,0) rectangle (\varR,-\vardi-\vardii);

	\draw (-\varR,-\vardi-\vardii) -- (\varR,-\vardi-\vardii);
	\draw (-\varR,0) -- (-\varr,0);
	\draw (\varr,0) -- (\varR,0);

	
	\draw [myorange, thick, fill=myorange!35] (-.1,0) ++(180:\varr) arc (180:360:\varr+0.1);
	\draw [myorange, thick, fill=myorange!20] (0,0) ++(180:\varr) arc (180:360:\varr);

	\foreach \angle in {-180, -150, -120, -90, -60, -30, 0}
	{
		\draw [->, myorange, thick] ({\varr*cos(\angle)}, {\varr*sin(\angle)}) -- ({0.75*\varr*cos(\angle)}, {0.75*\varr*sin(\angle)});
	}

\end{tikzpicture}
	\caption{Solidification interface and its velocity during the cooling stage (sectional view).}
	\label{fig:velocity}
\end{figure}

We begin by characterizing the velocity of a point $x(t)$ on some moving isothermal surface in $\Omega$; see \cref{fig:velocity}.
Since the temperature $\theta(x(t),t)$ is constant, we obtain
\begin{equation} \label{eq:dxt}
	\frac{\d}{\d t} \theta(x(t),t) 
	= 
	\grad \theta(x(t),t) \cdot x_t(t) + \theta_t(x(t),t) 
	= 
	0
	.
\end{equation}
The derivative $x_t(t)$ can be decomposed as
\begin{equation} \label{eq:xt_alpha_beta}
	x_t(t) 
	= 
	\alpha(x(t),t) \, \grad \theta(x(t),t) + \text{component perpendicular to } \grad \theta(x(t),t),
\end{equation}
where $\alpha(x(t),t)$ is a scalar function. 
Substituting \eqref{eq:xt_alpha_beta} into \eqref{eq:dxt}, we obtain 
\begin{equation}
	\alpha(x(t),t) 
	= 
	\frac{{}-{} \theta_t(x(t),t)}{\norm{\grad \theta(x(t),t)}^2}
	,
\end{equation}
where $\norm{\cdot}$ denotes the Euclidean norm.
Therefore, we can define the \emph{velocity of any isothermal surface} passing through the point~$x$ at time~$t$ as
\begin{equation*}
	v(x,t) 
	\coloneqq 
	\frac{{}-{} \theta_t(x,t)}{\norm{\grad \theta(x,t)}}
	.
\end{equation*}

While the melting pool expands, $v(x,t)$ takes negative values near the edge of the pool since $\theta_t > 0$ holds.
When the pool shrinks, $v(x,t)$ has positive values.
We are only interested in restricting positive velocities and only within the solidus--liquidus temperature corridor.
We therefore propose the following penalty term,
\begin{equation} \label{eq:J_velocity}
	J_\textup{velocity} 
	= 
	\frac{\beta_\textup{velocity}}{2} \int_{\Omega \times (0,T)} \max \paren[big]\{\}{v(x,t) - v_{\max}, \; 0}^2 \cdot \chi(\theta(x,t)) \d x \d t
	,
\end{equation}
where $v_{\max}$ is a predefined constant and the indicator function $\chi$ is defined as
\begin{equation*}
	\chi(\theta) 
	\coloneqq 
	\begin{cases}
		1 & \text{where}\ \text{solidus} \le \theta \le \text{liquidus}
		, 
		\\
		0 & \text{otherwise}
		.
	\end{cases}
\end{equation*}

\subsection{Completeness of Solidification}

To ensure that the \textbf{solidification stage is complete} at the given final time $T$, we penalize final temperatures $\theta(x,T)$ which are still above the solidus temperature by means of the following term,
\begin{equation} \label{eq:J_completeness}
	J_\textup{completeness} 
	=
	\frac{\beta_\textup{completeness}}{2} \int_{\Omega} \max \paren[big]\{\}{\theta(x, T) - \text{solidus}, \; 0}^2 \d x
	.
\end{equation}

\subsection{Energy Consumption Penalties}

The \textbf{consumption of energy} in the process is taken into account by means of the following standard quadratic control cost term,
\begin{equation}
	J_\textup{control} 
	=
	\frac{\beta_\textup{control}}{2} \norm{u}^2_{L^2(0,T)}
	.
\end{equation}
Indeed, an $L^1$-norm penalty would be a more meaningful model of energy consumption.
Such a term is known to induce sparsely supported controls, see for instance \cite{VossenMaurer:2006:1,Stadler:2009:1,CasasHerzogWachsmuth:2012:2}.
In the present application, however, optimal power profiles may then require the laser to be switched off and on again.
Technical limitations require a certain amount of time before the laser can be powered up again, which is not feasible due to the brevity of the usual process times~$T$ in single-spot welding.
Moreover, a waiting-time constraint would render the optimal control problem significantly more difficult.

\subsection{Optimal Control Problem Formulation}

For convenience, we summarize our single-spot welding optimal control problem as follows:

\begin{equation} \label{eq:J}
	\paren[auto].\}{%
		\begin{aligned}
			&
			\text{Find a control function $u \colon [0,T] \to \R$ which minimizes the objective}
			\\
			&
			J(u,\theta) 
			\coloneqq 
			J_\textup{penetration}(\theta) + J_\textup{velocity}(\theta) + J_\textup{completeness}(\theta) + J_\textup{control}(u)
			,
			\\
			&
			\text{where $\theta$ is the solution to the boundary value problem \eqref{eq:heat_eq}--\eqref{eq:heat_eq_bc}}
			\\
			&
			\text{and the control satisfies the constraints $0 \le u(t) \le 1$ on $[0,T]$}
			.
		\end{aligned}
	}
\end{equation}

\section{Discretization and Optimization Scheme}
\label{sec:discretization}

In this section we describe a discretization of problem \eqref{eq:J} as well as a projected gradient descent scheme for its numerical solution.
Since the discretization in space is based on a finite element approach, we begin with the notion of weak solution.
Notice that our definitions are informal since we do not aim to provide a thorough analysis of the forward system \eqref{eq:heat_eq}--\eqref{eq:heat_eq_bc} here.

\subsection{Weak Formulation}

As usual, the weak formulation is obtained by multiplying \eqref{eq:heat_eq} by a test function, integrating by parts, and plugging in the natural boundary conditions \eqref{eq:heat_eq_bc}.
Abbreviating
\begin{equation*}
	\Phi(\theta(x,t)) 
	\coloneqq 
	k \paren[big](){\theta(x,t)^4 - \theta_\textup{amb}^4} + h \paren[big](){\theta(x,t) - \theta_\textup{amb}}
	,
\end{equation*}
we thus arrive at the notion that a function $\theta \colon \Omega \times [0,T] \to \R$ is a \emph{weak solution} to the boundary value problem \eqref{eq:heat_eq}--\eqref{eq:heat_eq_bc} if it satisfies the initial condition \eqref{eq:heat_eq_ic} and the equality
\begin{multline} \label{eq:weak}
	\int_\Omega s(\theta(x,t)) \, \theta_t(x,t) \, v \d x \d t
	+
	\int_\Omega \grad \theta(x,t)^\transp \kappa(\theta(x,t)) \, \grad v \d x \d t 
	\\
	+
	\int_{\Gamma_1 \cup \Gamma_2 \cup \Gamma_4} \Phi(\theta(x,t))\, v \d S \d t -
	\int_{\Gamma_1} \eta \, \text{pd}_{\max} \, u(t) \, v \d S \d t 
	= 
	0
\end{multline}
holds for all functions $v \in C^\infty(\Omega)$ and for almost all $t \in (0,T)$.
Notice that $\d x$ denotes integration \wrt the volume measure and $\d S$ is \wrt the surface measure.
Recall that the thermal conductivity $\kappa(\theta)$ is a matrix due to different conductivities in radial and axial directions, see \cref{subsec:conductivity}.

\subsection{Reduction to the Radially Symmetric Case}

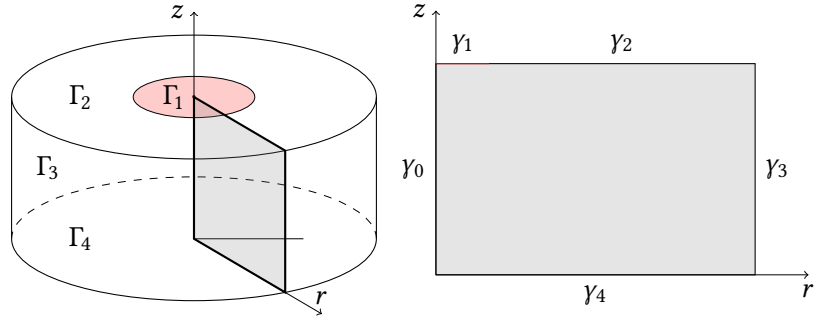
\begin{figure}
	\centering

\newcommand{\varR}{3.0}
\newcommand{\varr}{1.0}
\newcommand{\varZ}{2.5}

\tdplotsetmaincoords{70}{0}

\begin{tikzpicture}[tdplot_main_coords, scale=0.8]

	\begin{scope}[canvas is xy plane at z=0]
		\draw (\varR, 0) arc [radius=\varR, start angle=0, end angle=-180];
		\draw [very thin, dashed] (\varR, 0) arc [radius=\varR, start angle=0, end angle=180];
	\end{scope}

	\begin{scope}[canvas is xy plane at z=\varZ]
		\draw [fill] node{.} (0,0);
		\draw (0, 0) circle [radius=\varR];
		\draw [very thin, fill=red, fill opacity=0.2] (0, 0) circle [radius=\varr];
	\end{scope}

	\draw ( \varR, 0, 0) -- ( \varR, 0, \varZ);
	\draw (-\varR, 0, 0) -- (-\varR, 0, \varZ);

	\tdplotsetcoord{P1}{ .5*\varR}{90}{-60}
	\tdplotsetcoord{P2}{    \varR}{90}{-60}
	\tdplotsetcoord{P3}{1.4*\varR}{90}{-60}

	\draw [very thin] (0, 0, 0) -- (P2);
	\draw [->, very thin] (P2) -- (P3) node [above] {$r$};
	\draw [very thin] (0, 0, 0) -- (0, 0, \varZ);
	\draw [->, very thin] (0, 0, \varZ) -- (0, 0, 1.6*\varZ) node [left] {$z$};
	\draw [very thin] (0, 0, 0) -- (.6*\varR, 0, 0);

	\node[left] at (0, 0, \varZ) {$\Gamma_1$};
	\node[left] at (-.5*\varR, 0,  \varZ) {$\Gamma_2$};
	\node at (-.8*\varR, 0, .5*\varZ) {$\Gamma_3$};
	\node[left] at (-.5*\varR, 0,  0) {$\Gamma_4$};

	\coordinate (P2z) at ($ (P2) + (0, 0, \varZ)$);
	\draw [thick, fill=gray, fill opacity=0.2] (0, 0, 0) -- (P2) -- (P2z) -- (0, 0, \varZ) -- cycle;
	
\end{tikzpicture}
\begin{tikzpicture}[scale=0.7]
	\draw [->, very thin] (0, 0) -- (7, 0) node[below] {$r$};
	\draw [->, very thin] (0, 0) -- (0, 5) node[left] {$z$};

	\draw [fill=gray, fill opacity=0.2] (0, 0) -- (6, 0) -- (6, 4) -- (0, 4) -- cycle;

	\node[left] at (0, 2) {$\gamma_0$};
	\node[above] at (0.5, 4) {$\gamma_1$};
	\node[above] at (3.5, 4) {$\gamma_2$};
	\node[right] at (6, 2) {$\gamma_3$};
	\node[below] at (3, 0) {$\gamma_4$};

	\draw [red] (0, 4) -- (1, 4);
\end{tikzpicture}
	\caption{Draft: reduction to the radially symmetric case.}
	\label{fig:sec}
\end{figure}

Up to this moment the problem was considered in $\R^3$.
However, the power density of the laser beam is taken to be radially symmetric, and there is no heat transition in $\Omega$ in the angular direction, \ie $\partial\theta/\partial\varphi = 0$. 
This motivates us to reduce the computational complexity of the problem by reducing the domain $\Omega$ to its two-dimensional radial section $\omega$, see \cref{fig:sec}.

From now on, we replace $\theta(x,t)$ by $\theta(r,z,t)$; see \cref{fig:cylinder} for the coordinate axes.
Thus, \eqref{eq:weak} turns becomes
\begin{multline} \label{eq:weak_reduced}
	\int_\omega s(\theta(r,z,t)) \frac{\d \theta(r,z,t)}{\d t} v \, r \d r \d z \d t
	+
	\int_\omega \grad\theta(r,z,t)^\transp \kappa(\theta(r,z,t)) \grad v \, r \d r \d z \d t 
	\\
	+
	\int_{\gamma_1 \cup \gamma_2 \cup \gamma_4} \Phi(\theta(r,z,t)) \, v \, r \d s \d t -
	\int_{\gamma_1} \eta \, \text{pd}_{\max} \, u(t) \, v \, r \d s \d t 
	= 
	0
\end{multline}
for all $v \in C^\infty(\omega)$ and for almost all $t \in (0,T)$.
Notice that the gradient operator in equation~\eqref{eq:weak_reduced} must be used in its cylindrical form, \ie,
\begin{equation*}
	\grad \theta(r,z,\varphi) 
	=
	\frac{\partial \theta}{\partial r} e_r
	+
	\frac{\partial \theta}{\partial z} e_z
	+
	\frac{1}{r} \frac{\partial \theta}{\partial \varphi} e_{\varphi}
	.
\end{equation*}
However, as mentioned before, due to the radial symmetry of the heat distribution, the $e_{\varphi}$-component of $\theta$ vanishes. 
This feature is convenient for the numerical implementation, since the standard gradient operator (in Cartesian coordinate form) can be used.
In \eqref{eq:weak_reduced}, we now denote the surface measure of the two-dimensional cross-sectional domain~$\omega$ by $\d s$.
Notice that the integrals in \eqref{eq:weak_reduced} and in the following incur an extra factor~$r$ due to the coordinate transformation.

Similarly, two of the penalty terms in the objective in \eqref{eq:J} are affected by the transition to cylindrical coordinates.
Specifically, \eqref{eq:J_velocity} and \eqref{eq:J_completeness} now take the following forms:
\begin{align}
	\label{eq:J_velocity_r}
	J_\textup{velocity} 
	&
	=
	\frac{\beta_\textup{velocity}}{2} \int_{\omega \times (0,T)} \max \paren[big]\{\}{v(r,z,t) - v_{\max}, \; 0}^2 \cdot \chi(\theta(r,z,t)) \, r \d r \d z \d t
	, 
	\\
	\label{eq:J_completeness_r}
	J_\textup{completeness} 
	&
	=
	\frac{\beta_\textup{completeness}}{2} \int_{\omega} \max \paren[big]\{\}{\theta(r,z,T) - \text{solidus}, \; 0}^2 \, r \d r \d z
	.
\end{align}

\subsection{Discretization of the Forward Problem}

We now focus on discretizing the problem in time and space in order to solve it numerically. 
We combine a finite element method in space with a finite difference method in time. 
The numerical implementation is based on the \fenics computing platform; see \cite{LoggMardalWells:2012:1}.

Let $N_t$ be the number of equidistant time steps excluding the initial state, then we denote:
\begin{equation*}
	\begin{aligned}
		&
		\tau 
		\coloneqq 
		T/N_t
		, 
		\quad
		u_n 
		\coloneqq 
		u(n\tau)
		,
		\quad
		\theta_n(r,z) 
		\coloneqq 
		\theta(r,z,n\tau)
		, 
		\\
		&
		\theta_{n+\alpha}(r,z) 
		\coloneqq 
		\alpha \, \theta_{n+1}(r,z) + (1-\alpha) \, \theta_n(r,z)
		,
	\end{aligned}
\end{equation*}
where $\alpha \in [0,1]$ determines the degree of implicitness of the time scheme.

Within the time interval $(n\tau, n\tau+\tau]$, the coefficients and the operators of equation \eqref{eq:heat_eq} are discretized as follows:
\begin{equation*}
	\begin{aligned}
		s(\theta(r,z,t)) 
		&
		\coloneqq 
		s(\theta_n)
		, 
		&
		\kappa(\theta(r,z,t)) 
		&
		\coloneqq 
		\kappa(\theta_n)
		, 
		&
		\Phi(\theta(r,z,t)) 
		&
		\coloneqq 
		\Phi \paren[auto](){\theta_{n+\alpha}}
		, 
		\\
		\frac{\d\theta(r,z,t)}{\d t} 
		&
		\coloneqq 
		\frac{\theta_{n+1}-\theta_n}{\tau}
		, 
		&
		\grad(\theta(r,z,t)) 
		&
		\coloneqq 
		\grad \paren[auto](){\theta_{n+\alpha}}
		.
	\end{aligned}
\end{equation*}

For the discretization in space, we employ piecewise linear, globally continuous test and trial functions on a predefined mesh of $\omega$.
Now the discretized form of equation \eqref{eq:weak_reduced} reads as follows,
\begin{multline} \label{eq:discrete}
	\sum_{n=0}^{N_t-1} \int_{\omega}
	s(\theta_n) (\theta_{n+1}-\theta_n) \, v_n \, r \d r \d z
	+ \tau \sum_{n=0}^{N_t-1} \int_{\omega}
	\grad \theta_{n+\alpha}^\transp \, \kappa(\theta_n) \grad v_n \, r \d r \d z 
	\\
	+ \tau \sum_{n=0}^{N_t-1} \int_{\gamma_1 \cup \gamma_2}
	\Phi \paren[auto](){\theta_{n+\alpha}} v_n \, r \d s
	- \tau \sum_{n=0}^{N_t-1} \int_{\gamma_1}
	\eta \, \text{pd}_{\max} \, u_n \, v_n \, r \d s 
	= 
	0
	.
\end{multline}
In \eqref{eq:discrete} we set $\theta_0 \coloneqq \theta_\textup{amb}$.
We then solve \eqref{eq:discrete} time step by time step for the unknown coefficient vectors $\theta_1, \theta_2, \ldots, \theta_{N_t}$.

\subsection{Discretization of the Objective Functional}
\label{subsec:disc_objective}

To derive the discrete version of $J_\textup{penetration}$, we discretize the $L^p$-norm in \eqref{eq:J_penetration} according to
\begin{equation*}
	\norm{\theta(x_{\text{target}},\cdot)}_{L^p(0,T)}
	\approx
	\paren[Big](){\tau \sum_{n=1}^{N_t} \abs[big]{\theta_n(0,z_\textup{target})}^p}^{1/p}
	=
	\tau^{1/p}\ \norm[big]{ \{ \theta_n(0,z_\textup{target}) \}_{n=1}^{N_t} }_{l^p}
	.
\end{equation*}
In fact, the factor $\tau^{1/p}$, which tends to one as $p$ tends to infinity, can be compensated by adjusting the coefficient $\beta_\textup{penetration}$, so we implement the following discrete version of \eqref{eq:J_penetration}:
\begin{equation} \label{eq:J_penetration_discrete}
	J_\textup{penetration} 
	= 
	\frac{\beta_\textup{penetration}}{2}
	\paren[auto]\{\}{\paren[Big](){\sum_{n=1}^{N_t} \abs[big]{\theta_n(0,z_\textup{target})}^p}^{1/p} - \theta_\textup{target}}^2
	.
\end{equation}
A detailed discussion on the choice of $p$ is given in \cref{subsec:p-norm_discussion}.

The velocity of an isothermal surface (in fact an isothermal line after dimension reduction) can be approximated as
\begin{equation*}
	v(\theta_n, \theta_{n+1}) 
	= 
	\frac{{}-{}(\theta_{n+1}-\theta_n)}{\tau \, \norm{\grad \theta_{n+\alpha}}}
\end{equation*}
and hence $J_\textup{velocity}$ takes the following form:
\begin{equation}
	J_\textup{velocity} 
	=
	\frac{\beta_\textup{velocity}}{2} \tau \sum_{n=0}^{N_t-1} \int_{\omega} \max \paren[big]\{\}{v(\theta_n, \theta_{n+1}) - v_{\max}, \; 0}^2 \cdot
	\chi(\theta_n, \theta_{n+1})
	\, r \d r \d z
\end{equation}
where the discretized indicator function $\chi$ is defined as
\begin{equation*}
	\chi(\theta_n, \theta_{n+1})
	\coloneqq
	\begin{cases}
		1 & \text{where}\ \text{solidus} \le \theta_n\ \text{and}\ \theta_{n+1} < \text{liquidus}
		,
		\\
		0 & \text{otherwise}
		.
	\end{cases}
\end{equation*}

The remaining penalty terms $J_\textup{completeness}$ and $J_\textup{control}$ are discretized according to
\begin{align}
	J_\textup{completeness} 
	&
	=
	\frac{\beta_\textup{completeness}}{2} \int_{\omega} \max \paren[big]\{\}{\theta_{N_t} - \text{solidus}, \; 0}^2 \, r \d r \d z
	, 
	\\
	J_\textup{control} 
	&
	=
	\frac{\beta_\textup{control}}{2} \tau \sum_{n=0}^{N_t-1} u_n^2
	.
\end{align}

\subsection{Evaluation of the Gradient}

In this section we briefly describe the evaluation of the gradient by means of the discrete adjoint state $p = [p_0, p_1, \ldots, p_{N_t-1}]$.
To this end, we introduce the Lagrangian 
\begin{multline} \label{eq:lagrange}
	\cL(\theta,u,p) \coloneqq
	J(\theta,u)
	+
	\sum_{n=0}^{N_t-1}\int_{\Omega} s(\theta_n) (\theta_{n+1}-\theta_n) \, p_n \, \d x
	+ 
	\tau \sum_{n=0}^{N_t-1}\int_{\Omega} \grad \theta_{n+\alpha}^\transp \, \kappa(\theta_n) \grad p_n \d x
	\\
	+ 
	\tau \sum_{n=0}^{N_t-1}\int_{\gamma_1 \cup \gamma_2} \Phi \paren[auto](){\theta_{n+\alpha}} p_n \d s
	- 
	\tau \sum_{n=0}^{N_t-1}\int_{\gamma_1} \eta \, \text{pd}_{\max} \, u_n \, p_n \d s.
\end{multline}
The sequence of linear systems governing the discrete adjoint state is obtained from $\partial \cL(\theta,u,p)/\partial \theta_n = 0$.
We do not provide the explicit formula for the adjoint equation here since in the code we derive it using \fenics' built-in automatic differentiation capabilities.
The only manual differentiation required is for the penalty term $J_\textup{penetration}$ in \eqref{eq:J_penetration_discrete}, since in contrast to the other terms, it cannot be split into a sum over the time steps.
We added the contributions coming from this term manually to the adjoint state's right hand side.
One can find more details in \cite[\texttt{optipuls.core}]{optipuls_github}.

Finally, we differentiate $\cL(\theta,u,p)$ with respect to $u = [u_0, u_1, \ldots, u_{N_t-1}]$ in the direction $\delta u$ to obtain
\begin{equation*}
	\frac{\partial \cL(\theta,u,p)}{\partial u} \, \delta u
	=
	\tau \sum_{n=0}^{N_t-1}
	\paren[auto][]{\beta_\textup{control} \, u_n - \int_{\gamma_1} \eta \, \text{pd}_{\max} \, p_n \d s} \, \delta u
	.
\end{equation*}
Consequently,
\begin{equation*} \label{eq:gradient}
	\grad_u \cL(\theta, u, p) 
	= 
	\beta_\textup{control}\, u - \int_{\gamma_1} \eta \, \text{pd}_{\max} p \d s
\end{equation*}
holds.

\subsection{Projected Gradient Descent Scheme}

To find the optimal control for the discretized counterpart of \eqref{eq:J}, we apply a \emph{projected gradient descent} scheme with line search; see, \eg, \cite{GafniBertsekas:1984:1,CalamaiMore:1987:1} or \cite[Chapter~5.8.2]{GeigerKanzow:2002:1}.
To this end, we denote by $j(u) = J(u,\theta)$ the reduced objective, which depends only on the values $u = [u_0, u_1, \ldots, u_{N_t-1}]$ of the control since the solution $\theta = [\theta_1, \theta_2, \ldots, \theta_{N_t}]$ to the forward system \eqref{eq:discrete} has been inserted.
Since this procedure is well known, we present only a short general outline in \cref{alg:projected_gradient_descent}.
The norm in which the size of the gradient is evaluated is the norm represented by $\tau$ times the identity matrix.
More details can be found in the implementation at \cite[\texttt{optipuls.optimization}]{optipuls_github}.
 
\begin{algorithm2e}[H]
	\caption{Projected gradient descent scheme.}
	\label{alg:projected_gradient_descent}
	\DontPrintSemicolon
	\KwIn{$u_\textup{initial} \in \R^{N_t}$}
	\KwOut{$u_\textup{optimized} \in \R^{N_t}$}
	$u_\textup{current} \gets u_\textup{initial}$\;
	\While{stopping criterion is not satisfied}{
		solve the forward system \eqref{eq:discrete} for $\theta$, given $u_\textup{current}$\;
		solve the adjoint system for $p$, given $\theta$ and $u_\textup{current}$\;
		evaluate the gradient of the reduced objective $\grad_u j(u_\textup{current})$ from \eqref{eq:gradient}\;
		\Repeat{$J(u_\textup{trial}) \ge J(u_\textup{current}) - \sigma \, \alpha \, \norm{\grad_u j (u_\textup{current})}^2$}{
			perform a line search for the step size $\alpha$\;
			$u_\textup{trial} \gets \cP_{[0,1]} \paren[auto](){u_\textup{current} - \alpha \grad_u j (u_\textup{current})}$\;
		}
		$u_\textup{current} \gets u_\textup{trial}$\;
	}
	\Return{$u_\textup{optimized} \gets u_\textup{current}$}
\end{algorithm2e}

The stopping criterion was considered satisfied as soon as any of the following conditions were met.
\begin{equation*}
	\begin{aligned}
		\norm{\cP_A \grad_u j (u_\textup{current})}
		&
		< 
		\text{tolerance}_{\grad}
		, 
		\\
		\norm{u_\textup{trial} - u_\textup{current}}
		&
		< 
		\text{tolerance}_\text{control}
		, 
		\\
		1 - J(u_\textup{trial}) / J(u_\textup{current})
		&
		< 
		\text{tolerance}_\text{descent rate}
		, 
		\\
		\text{iteration no.}
		&
		> 
		M
		.
	\end{aligned}
\end{equation*}
Here $\cP_A$ is the point-wise projection onto the tangent cone to the feasible set in $L^2(0, T; [0,1])$ at $u_\textup{current}$, \ie,
\begin{equation*}
	A
	= 
	\begin{cases}
		(-\infty, 0]
		&
		\text{where }
		u_\textup{current} = 0
		,
		\\
		(-\infty, \infty)	
		& 
		\text{where }
		0 < u_\textup{current} < 1
		,
		\\
		[0, \infty)
		& 
		\text{where }
		u_\textup{current} = 1
		.
	\end{cases}
\end{equation*}

\section{Numerical Results}
\label{sec:numericals}

In this section we present some optimized laser pulses, \ie, numerical solutions to the discretized counterpart of the single-spot welding optimal control problem \eqref{eq:J}.
We emphasize that all numerical results presented in this paper are fully reproducible and hence can be verified by the reader; see \cite{optcontrol_github} for further instructions.

The common problem parameters shared by the numerical experiments presented in this section are provided in \cref{tab:parameters} (see also \cite[\texttt{env/problem.py}]{optcontrol_github}).
The parameters are describing the EN~AW~6082-T6 aluminum alloy and were kindly provided by the Department of Production Technology, TU~Ilmenau, Germany.

Notice that the target temperature $\theta_\text{target}$ in \cref{tab:parameters} is intentionally set higher than the desired maximal temperature at the target point ($\text{liquidus} = \SI{923}{\K}$), see \cref{subsec:p-norm_discussion} for further details.

\begin{table}[ht]
	\begin{tabular}{rl}
		\textbf{space domain} \\
		radius of the cylinder $\Omega$ & $\SI{2.5}{\mm}$ \\
		height of the cylinder $\Omega$ & $\SI{0.5}{\mm}$ \\
		radius of the laser beam        & $\SI{0.2}{\mm}$\\[0.5em]

		\textbf{equation and boundary conditions} \\
		ambience temperature & $\theta_\text{amb} = \SI{295}{\K}$ \\
		convective heat transfer coefficient & $h = 20\, \si{\W\per\m^2}$\\
		radiative heat transfer coefficient & $k = 2.26 \cdot 10^{-9}\, \si{\W\per\m^2\K^4}$\\[0.5em]

		\textbf{objective functional} \\
		welding penetration penalty coefficient & $\beta_\text{penetration} = 10^{-2}$ \\
		solidification velocity penalty coefficient & $\beta_\text{velocity} = 1.5 \cdot 10^{-1}$ \\
		welding completeness penalty coefficient & $\beta_\text{completeness} = 10^{-12}$ \\
		energy consumption penalty coefficient & $\beta_\text{control} = 10^{2}$ \\
		target point & $z_\text{target} = \SI{0.375}{\mm}$ (welding depth $\SI{0.125}{\mm}$) \\
		target maximal temperature at the target point & $\theta_\text{target} = \SI{1048}{\K}$ \\
		$p$-norm in time domain& $p = 20$\\[0.5em]

		\textbf{material properties} \\
		solidus point & $\SI{858}{\K}$ \\
		liquidus point & $\SI{923}{\K}$ \\
		enthalpy of fusion & $\SI{397000}{\J\per\kg}$ \\
		coefficients $s(\theta)$ and $\kappa(\theta)$ & discussed in \cref{sec:modelling}, see \cref{fig:coef} \\
	\end{tabular}
	\caption{Parameters of the numerical experiments}
	\label{tab:parameters}
\end{table}

\subsection{Conventional and Linear Rampdown Pulse Shapes}

Conventional pulsed laser welding strategies use a rectangular laser pulse shape, \ie the laser is working full power for a short time and is switched off immediately after. 
Unfortunately, this simple strategy often leads to hot cracking when applied to aluminum alloys.
A so-called linear rampdown pulse shape, \ie when the laser power is decreasing linearly after a short period of working full power, has shown its potential to obtain a crack-free welding of aluminum alloys; see \cite{ZhangWeckmanZhou:2008:1,JiaZhangYuShiLiuWuYeWangTian:2021:1}.
However, rampdown pulses are not likely to be optimal with respect to any of the criteria established in \cref{sec:optimal_control_problem}.

\begin{figure} 
	\centering
	\includegraphics{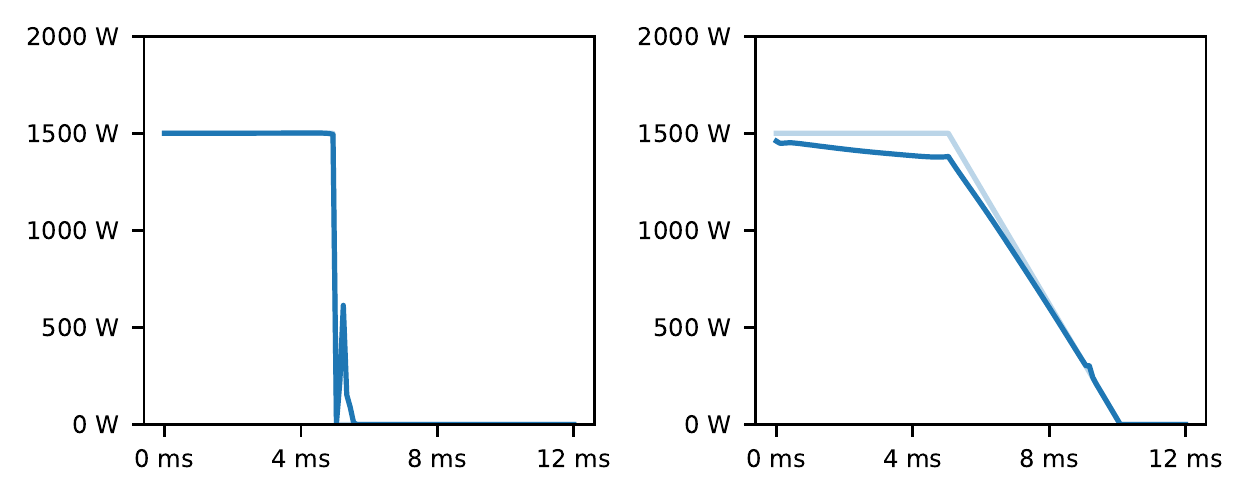}
	\caption{Solutions to the optimal control problem with conventional (left) and linear rampdown (right) pulse shapes taken as initial guesses.}
	\label{fig:rampdown}
\end{figure}

In view of this we first consider the conventional ($\SI{5}{\ms}$ of $\SI{1500}{\W}$) and the linear rampdown ($\SI{5}{\ms}$ of $\SI{1500}{\W}$, $\SI{5}{\ms}$ of rampdown) pulse shapes as the initial guess $u_\textup{initial}$ for the optimizer. In both of these experiments the maximal laser power $P_\text{YAG}$ is limited by $\SI{2000}{\W}$ and the total time $T$ is limited by $\SI{12}{\ms}$.
\Cref{fig:rampdown} demonstrates the corresponding solutions to the optimal control problem obtained with \cref{alg:projected_gradient_descent}.
The numerical reports on the corresponding simulations are presented in \cref{tab:rampdown}. Notice that an imperfect match between the evaluated welding depth and the $J_\text{penetration}$ penalty is due to the $p$-norm approximation to the sup-norm of $\theta(x_\text{target}, \cdot)$, see~\eqref{eq:p-norm_ineq}.
We observe that in both cases, we obtain apparently locally optimal pulse shapes which differ very little from their respective initial guesses.

\begin{table} 
	\centering
	\begin{tabular}{lcccccc}
\toprule
 pulse shape       &  welding depth  &  $J_\text{penetration}$  &  $J_\text{velocity}$  &  $J_\text{completeness}$  &  $J_\text{control}$  &  $J_\text{ total}$  \\
\midrule
 conventional      &     0.11875     &          0.9758          &       278.9010        &          0.0000           &        0.1406        &      280.0173       \\
 linear rampdown   &     0.11875     &          2.4323          &        0.0055         &          0.0000           &        0.1889        &       2.6267        \\
 conventional *    &     0.11875     &          0.8611          &        27.2947        &          0.0000           &        0.1413        &       28.2972       \\
 linear rampdown * &     0.11250     &          0.0001          &        0.0000         &          0.0000           &        0.1674        &       0.1675        \\
\bottomrule
\end{tabular}

	\caption{Numerical report on simulations with the conventional and the linear rampdown pulse shapes, and locally optimal pulse shapes * starting from corresponding initial guesses.}
	\label{tab:rampdown}
\end{table}

\subsection{Optimizations from Zero Initial Guess}
\label{subsec:optimization_zeroguess}

\begin{figure} 
	\centering
	\includegraphics{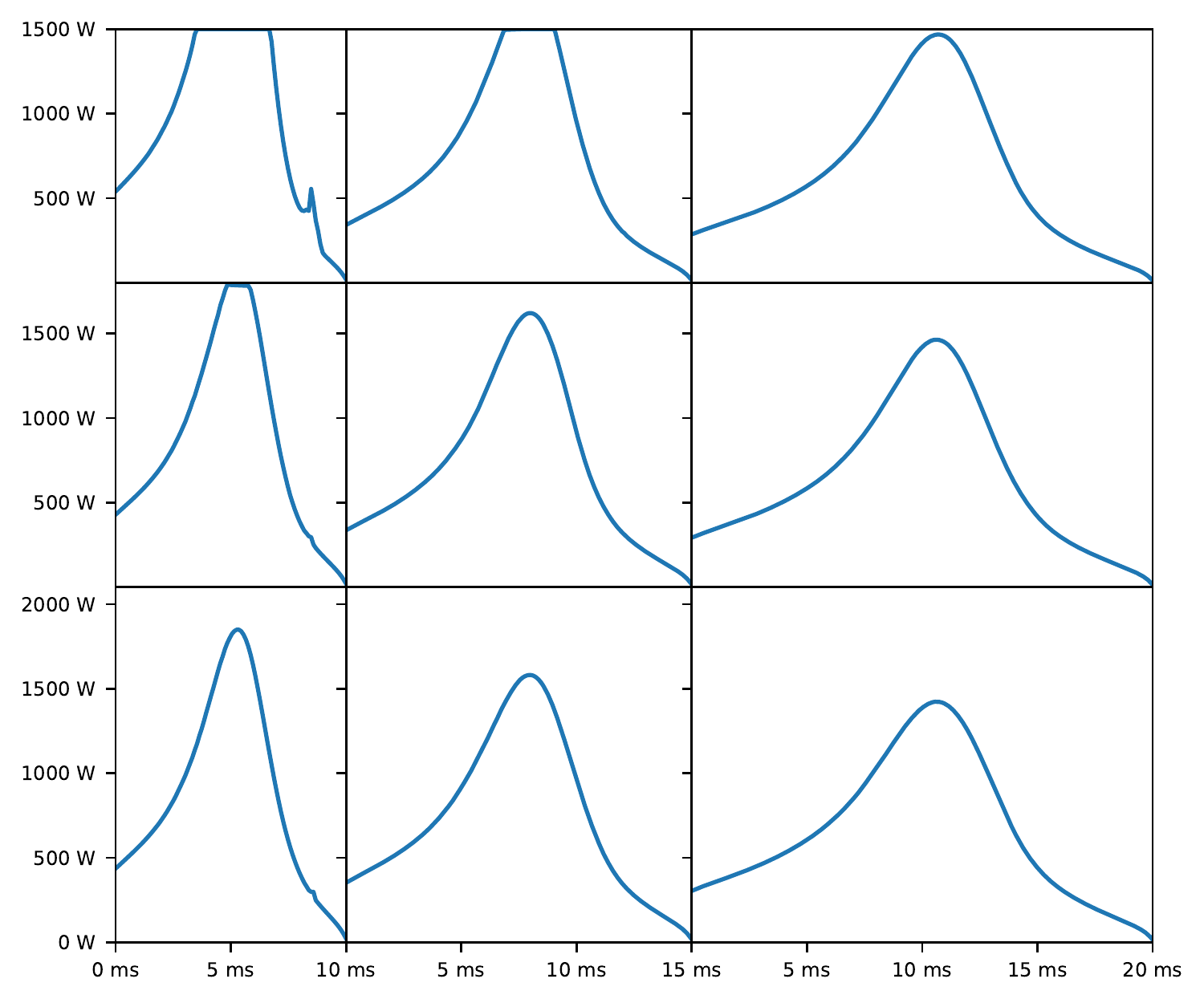}
	\caption{Solutions to the optimal control problem with zero initial guess with the maximal laser power varying vertically and the maximal welding time varying horizontally.}
	\label{fig:zeroguess}
\end{figure}

In the search to obtain better pulse shapes, we now begin with the trivial initial guess $u_\textup{initial} \equiv 0$, \ie, no power radiated by the laser. 
With this initial guess, the target temperature is clearly not reached and the term $J_\textup{penetration}$, see \eqref{eq:J_penetration_discrete}, drives the pulse shape away from its initial value.
\Cref{fig:zeroguess} shows the corresponding solutions to the optimal control problem with variable maximal laser power $P_\text{YAG}$ and maximal time $T$ until a full solidification.
The corresponding numerical reports are presented in \cref{tab:zeroguess}.

To give the reader some idea on the performance of the optimizer, the runtime of the gradient descent procedure was measured for the problem corresponding to $P_\text{YAG} = \SI{1800}{\W}$, $T = \SI{15}{\ms}$.
The descent terminated after 15~iterations by the \enquote{descent rate} stopping criterion and it took $\SI{175}{\s}$.
The temperature state had 4837 spatial degrees of freedom.
The computations were carried out on a workstation with an AMD Ryzen~9 5950X CPU.

\begin{table} 
	\centering
	\begin{tabular}{cccccccc}
\toprule
  $P_\text{YAG}$  &  $T$  &  welding depth  &  $J_\text{penetration}$  &  $J_\text{velocity}$  &  $J_\text{completeness}$  &  $J_\text{control}$  &  $J_\text{total}$  \\
\midrule
       1500       & 0.010 &     0.11875     &          0.0243          &        0.0958         &          0.0000           &        0.2522        &       0.3724       \\
       1500       & 0.015 &     0.11875     &          0.0002          &        0.0000         &          0.0000           &        0.2530        &       0.2532       \\
       1500       & 0.020 &     0.11250     &          0.0000          &        0.0000         &          0.0000           &        0.2640        &       0.2640       \\
       1800       & 0.010 &     0.12500     &          0.0000          &        0.0000         &          0.0000           &        0.1639        &       0.1640       \\
       1800       & 0.015 &     0.11875     &          0.0001          &        0.0000         &          0.0000           &        0.1741        &       0.1741       \\
       1800       & 0.020 &     0.11250     &          0.0000          &        0.0000         &          0.0000           &        0.1846        &       0.1846       \\
       2100       & 0.010 &     0.12500     &          0.0001          &        0.0000         &          0.0000           &        0.1206        &       0.1207       \\
       2100       & 0.015 &     0.11875     &          0.0000          &        0.0000         &          0.0000           &        0.1302        &       0.1302       \\
       2100       & 0.020 &     0.11250     &          0.0001          &        0.0000         &          0.0000           &        0.1375        &       0.1375       \\
\bottomrule
\end{tabular}

	\caption{Numerical report on the series of optimizations with zero initial guess.}
	\label{tab:zeroguess}
\end{table}

\subsection{Impact of the p-norm approximation to the sup-norm on the temperature control}
\label{subsec:p-norm_discussion}

The choice of the value of $p$ in the $p$-norm approximation to the sup-norm has significant impact on the accuracy within which the maximal temperature at the target point $x_{\text{target}}$ can be controlled via penalty term~\eqref{eq:J_penetration}. 

\begin{figure}
	\centering
	\includegraphics{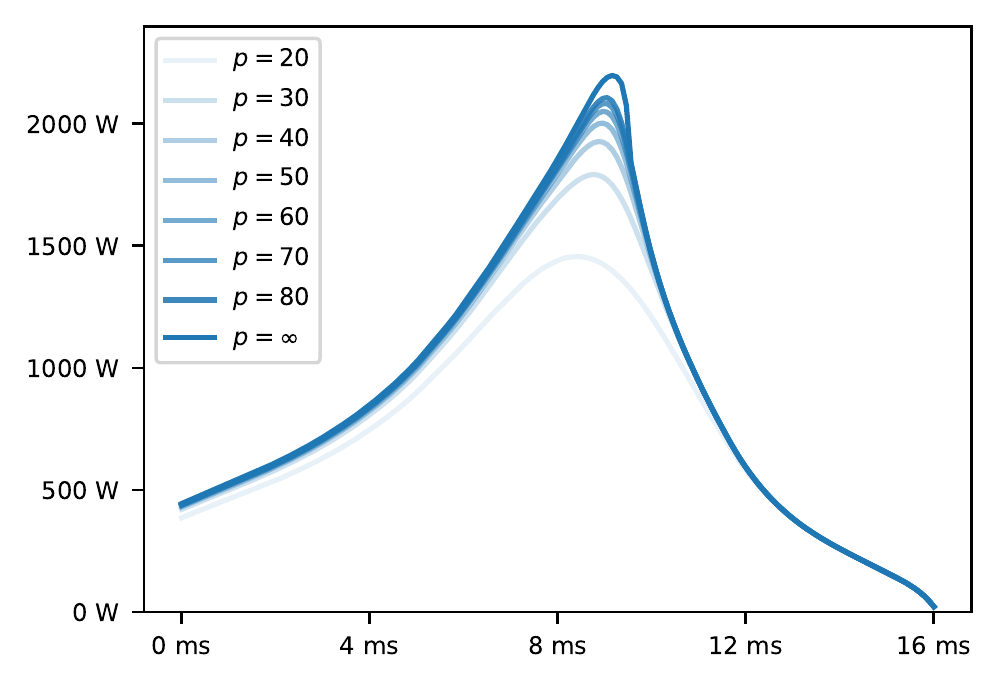}
	\caption{Solutions to a sequence of optimal control problems with variable parameter $p$. Zero initial guess is taken for the smallest value of $p$ and each subsequent problem utilizes previously computed optimal control as the initial guess.}
	\label{fig:p-norm}
\end{figure}

\cref{fig:p-norm} shows solutions to a sequence of optimal control problems employing successively increasing values of $p$ and \cref{tab:p-norm} shows the actual maximal temperature reached at the target point.
The zero initial guess is taken for the smallest value of $p$ and each subsequent problem utilizes the previously computed optimal control as the initial guess.

Although higher values of $p$ bring more accurate control of the welding penetration, they lead to a side effect which is rather undesirable for the practical application: solutions to the corresponding optimal control problems tend to demonstrate faster growth and higher peaks comparing to those evaluated for smaller values of $p$.
Due to higher laser beam power density such power profiles can form a so-called \emph{key hole} structure and therefore trigger a transition from the heat conduction welding to the deep penetration welding, see \cite{ZhouTsai:2013:1}. While deep penetration laser welding has its own applications, such transitions are strictly avoided in the current study.

Moreover, higher values of $p$ used in optimizations from zero initial guess described in \cref{subsec:optimization_zeroguess} lead to solutions showing higher solidification velocity rates and corresponding penalties. Therefore our strategy is to compensate smaller values of $p$ by setting the target temperature $\theta_\text{target}$ higher than the actual desired maximal temperature at the target point. As one can see from the following inequality
\begin{equation} \label{eq:p-norm_ineq}
	\norm[big]{ \{ \theta_n(0, z_\text{target}) \}_{n=1}^{N_t} }_{l^p}
	\le
	N_t^{\frac{1}{p}}\ \norm[big]{ \{ \theta_n(0, z_\text{target}) \}_{n=1}^{N_t} }_{l^{\infty}},
\end{equation}
it is enough to set $N_t^{1/p} \cdot \text{liquidus}$ as the target temperature in~\eqref{eq:J_penetration_discrete} to ensure that the liquidus temperature was reached at the target point, however the exact value was chosen by trial and error.


\begin{table}
	\centering
	\begin{tabular}{lrrrrrrrr}
\toprule
 $p$             &  20 &    30 &    40 &     50 &     60 &     70 &     80 & $\infty$ \\
 $\theta_{\max}$ & 879 & 955.2 & 986.6 & 1003.3 & 1013.5 & 1020.4 & 1025.2 &   1048.0 \\
\bottomrule
\end{tabular}

	\caption{Maximal temperature at the target point depending on variable parameter $p$.}
	\label{tab:p-norm}
\end{table}

\subsection{Discussion of the Obtained Numerical Results}

One can see from \cref{tab:rampdown} that the conventional pulse shape is far from being optimal since it leads to enormous solidification velocity penalties. 
Even its local optimization does not achieve a significant improvement.
On the other hand, the linear rampdown pulse shape achieves a reduction of the solidification velocity and its local optimization is even able to keep this velocity within the permitted limit~$v_{\max}$ so that $J_\textup{velocity}$ is zero.
A successful crack-free welding of aluminum alloys using the linear rampdown pulse shape was confirmed experimentally, see \cite{ZhangWeckmanZhou:2008:1,JiaZhangYuShiLiuWuYeWangTian:2021:1}.

However, the results of optimizations with zero initial guess in \cref{tab:zeroguess} show that further optimization is still possible and the optimal pulse shapes are quite non-trivial to guess by trial and error. 
One can see from \cref{tab:zeroguess} and \cref{fig:zeroguess} that reasonably small penalty values can be obtained only if the optimizer has enough room to adjust the pulse shape in time and power dimensions.

In the pulse shapes limited by $\SI{10}{\ms}$, the small swing-ups close to the end appear as a result of compensation of the too high solidification velocity. 
Despite of the fact that such pulse shapes are local optimizers for the discrete version of \eqref{eq:J}, they do not seem reasonable for the practical application. 
For the optimization problem under consideration with the setup as in \cref{tab:parameters}, the most promising optimal pulse shape would be the one obtained with a maximal laser power of $\SI{2100}{\W}$ and total time of $\SI{15}{\ms}$.

Preliminary laboratory experiments carried by the Department of Production Technology, TU~Ilmenau, confirmed that the pulse shapes optimized from zero initial guess can indeed be used to produce crack-free welds. 
These results will be presented in a separate paper.

\printbibliography

\end{document}